\documentclass[a4paper, reqno, english,11pt]{article}

\usepackage{amsfonts, amsthm, amsmath, amssymb}

\RequirePackage{mathrsfs} \let\mathcal\mathscr

\numberwithin{equation}{section}

\renewcommand{\phi}{\varphi}

\newtheorem{theorem}{Theorem}
\newtheorem*{thm*}{Theorem}

\renewcommand{\le}{\leqslant}

\renewcommand{\ge}{\geqslant}

\theoremstyle{definition}

\newcommand{\dif}{\mathrm{d}}

\begin{document}

\title{Exponential sums with coefficients of certain Dirichlet series}
\author{Stephan Baier}
\maketitle

\begin{abstract}
Under the generalized Lindel\"of Hypothesis in the $t$- and $q$-aspects, we bound exponential sums with coefficients of Dirichlet series belonging to a certain class. 
We use these estimates to establish a conditional result on squares of Hecke eigenvalues at Piatetski-Shapiro primes.
\end{abstract}



\section{General assumptions}
In this paper, we derive a conditional estimate for exponential sums of the form
$$
\sum\limits_{n\sim N} a_ne(f(n)),
$$
where $a_n$ is the $n$-th coefficient of Dirichlet series $F(s)$ whose twists with Dirichlet characters satisfy the generalized Lindel\"of Hypothesis in the $t$- and $q$-aspects, and $f(x)$ is a function 
having certain properties. As an application, we
consider squares of Hecke eigenvalues at Piatetski-Shapiro primes. In the following, we state the required conditions on $F(s)$ and $f(x)$.\\ 

{\bf Conditions on the $L$-function:}\medskip\\
We assume that 
$$
F(s)=\sum_{n=1}^\infty a_nn^{-s}
$$
is a Dirichlet series, absolutely convergent for $\Re s>1$, which satisfies the following conditions a), b) and c).\medskip\\
a) $F(s)$ lies in the extended Selberg class of Dirichlet series which don't necessarily possess a functional equation, i.e. $F(s)$ has the following properties.\medskip

(i) ({\it Analiticity}) There exists some $m\in \mathbb{N}$ such that $(s-1)^mF(s)$ extends to an entire function of finite order.

(ii) ({\it Ramanujan conjecture}) $a_1=1$ and $a_n \ll_\varepsilon n^\varepsilon$ for any $\varepsilon>0$.

(iii) ({\it Euler product}) For $\Re s>1$, the function $F(s)$ can be written as a product over primes in the form
$$
F(s)=\prod_p F_p(s),
$$
where  $\log F_p(s)$ is a Dirichlet series of the form
$$
\log F_p(s)=\sum_{n=0}^\infty b_{p^k}p^{-ks}
$$
with complex coefficients $b_{p^k}$ satisfying
$$
b_{p^k}=O(p^{k\theta})
$$
for some $\theta<1/2$.\\ 
b) For any Dirichlet character $\chi$ define 
$$
F(s,\chi):=\sum\limits_{n=1}^{\infty} a_n \chi(n) n^{-s} \quad \mbox{for } \Re s>1.
$$
Then  $(s-1)^mF(s,\chi)$ extends to an entire function again. \\
c) The family of functions $F(s,\chi)$ satisfies the Lindel\"of Hypothesis in the $t$- and $q$-aspects, i.e. 
\begin{equation} \label{LH}
F\left(\frac{1}{2}+it,\chi\right) \ll |tq|^{\varepsilon} \quad \mbox{for all $|t|\ge 1$, $q\in \mathbb{N}$ and characters 
$\chi$ mod $q$.}  
\end{equation}\\

{\bf Conditions on $f$:}\medskip\\
We assume that $f:[1,\infty)\rightarrow \mathbb{R}$ satisfies the following conditions a)-f).\medskip\\
a) $f$ is three times continuously differentiable.\\
b) $f$ is monotonically increasing.\\
c) $f(x)\asymp f(2x)$ for all $x\ge 1$.\\
d) $f^{(k)}(x) \asymp f(x)/x^k \mbox{ for all } x\ge 1 \mbox{ and } k=1,2,3$.\\
e) $f'(x)+xf''(x)\asymp f(x)/x \mbox{ for all } x\ge 1$.\\
f) $2f''(x)+xf'''(x)\asymp f(x)/x^2 \mbox{ for all } x\ge 1$.\medskip\\

\section{Results}
Our first result is the following.

\begin{theorem} \label{expsum} Fix $\eta>0$. Suppose that $1\le N<N'\le 2N$ and
\begin{equation} \label{anotherfcond}
N^{1/2+\eta}\le f(N)\le N^{3/2-\eta}.
\end{equation}
Then, under the conditions in section 1, we have 
$$
\sum\limits_{N<n\le N'} a_n e(f(n))\ll_{f,\eta,\varepsilon} N^{19/22+\varepsilon}f(N)^{1/11}.
$$
\end{theorem}

We note that the above bounds are non-trivial if $\varepsilon<\eta/11$ and $N$ is large enough.

With applications in mind, we also prove the following modification of Theorem \ref{expsum}.
\begin{theorem} \label{expsummod} Let $m\in \mathbb{N}$. Then, under the conditions in Theorem \ref{expsum} and section 1, we have 
$$
\sum\limits_{\substack{n\sim N\\ (n,m)=1}} \mu^2(n)a_n e(f(n))\ll_{f,\eta,\varepsilon} m^{\varepsilon}N^{19/22+\varepsilon}f(N)^{1/11}.
$$
\end{theorem}

Let now $G$ be a Hecke eigenform of weight $\kappa$ for the full modular group SL$_2(\mathbb{Z})$. By $\lambda(n)$ we denote the normalized $n$-th Fourier coefficient of 
$G$, {\it i.e.}
$$
G(z)=\sum\limits_{n=1}^{\infty} \lambda(n)n^{(\kappa-1)/2}e(nz) \quad \mbox{for } \Im z>0, \mbox{ and } \lambda(1)=1.  
$$
These Hecke eigenvalues satisfy the multiplicative property
\begin{equation} \label{multprop}
\lambda(mn)=\sum\limits_{d|\mbox{\scriptsize gcd}(m,n)} \mu(d)\lambda\left(\frac{m}{d}\right)\lambda\left(\frac{n}{d}\right) \quad \mbox{for all } m,n\in\mathbb{N}
\end{equation}
and the Ramanujan conjecture $\lambda(n)\ll_{\varepsilon} n^{\varepsilon}$.

Let $L(\mbox{Sym}^2\ G,s)$ be the symmetric square $L$-function for $G$, defined by
$$
L(\mbox{Sym}^2\ G,s)=\zeta(2s)\sum\limits_{n=1}^{\infty} \lambda\left(n^2\right)n^{-s} \quad \mbox{for } \Re s>1.
$$
We note that by multiplying out the right-hand side, we get
\begin{equation} \label{theas}
L\left(\mbox{Sym}^2\ G,s\right)=\sum\limits_{n=1}^{\infty} a_nn^{-s} \quad \mbox{for } \Re s>1,
\end{equation}
where the coefficients of the Dirichlet series on the right-hand side satisfy $a_n=\lambda\left(n^2\right)$ for any squarefree $n$. 
Moreover, it is well-known that $L(\mbox{Sym}^2\ G,s)$ lies in the Selberg class and hence satisfies condition a) in section 1.

More generally, for any Dirichlet character $\chi$ let $L(\mbox{Sym}^2\ G\otimes \chi,s)$ be 
the symmetric square $L$-function for $G$ twisted with $\chi$, defined by
$$
L(\mbox{Sym}^2\ G \otimes \chi,s)=\sum\limits_{n=1}^{\infty} \chi(n)a_nn^{-s}=
L\left(2s,\chi^2\right)\sum\limits_{n=1}^{\infty} \chi(n) \lambda\left(n^2\right)n^{-s} \quad \mbox{for } \Re s>1.
$$
As a consequence of the work of Shimura \cite{Shim}, 
$L\left(\mbox{Sym}^2\ G \otimes \chi,s\right)$ extends analytically to the whole complex plane and hence satisfies condition b) 
in section 1. If $\chi$ is primitive, then $L\left(\mbox{Sym}^2\ G \otimes \chi,s\right)$ even lies in the Selberg class.

The Lindel\"of Hypothesis in the $t$- and $q$-aspects for the family of $L$-functions $L\left(\mbox{Sym}^2\ G \otimes \chi,s\right)$,
with $G$ fixed, asserts that
\begin{equation} \label{LHLG}
L\left(\mbox{Sym}^2\ G \otimes \chi,\frac{1}{2}+it\right)\ll (tq)^{\varepsilon} \quad \mbox{for all $|t|\ge 1$, $q\in \mathbb{N}$ and characters 
$\chi$ mod $q$.}  
\end{equation}
We note that it can be deduced from Theorem 1 in \cite{ConGho} that \eqref{LHLG} holds if $L\left(\mbox{Sym}^2\ G \otimes \chi,s\right)$
satisfies the Riemann Hypothesis for all primitive characters $\chi$.

In \cite{baierzhao},
we bounded the average of $\lambda(p)$ at Piatetski-Shapiro primes, i.e. primes of the form $p=\left[n^c\right]$ with
$n\in \mathbb{N}$ and $c>1$ fixed. The $c$-range for which we obtained a non-trivial result was $1<c<8/7$. In this range, we proved that
\begin{equation} \label{piatbound}
\sum\limits_{\substack{n\le N\\ \left[n^c\right]\in \mathbb{P}}} \lambda\left(\left[n^c\right]\right) \ll N\exp\left(-C\sqrt{\log N}\right),
\end{equation}
where here as in the sequel, $\mathbb{P}$ is the set of primes.
We posed the question if also an asymptotic estimate for the average of the {\it squares} of these Hecke eigenvalues at 
Piatetski-Shapiro can be established.  Employing Theorem \ref{expsummod}, we shall prove the following conditional result.

\begin{theorem} \label{heckesquares} Let $1<c< 25/24$ be fixed and $\mathbb{P}$ be the set of primes. Assume that \eqref{LHLG} holds.
Then we have
$$
\sum\limits_{\substack{n\le N\\ \left[n^c\right]\in \mathbb{P}}} \lambda\left(\left[n^c\right]\right)^2\sim \frac{N}{c\log N} \quad
\mbox{as } N\rightarrow \infty.
$$ 
\end{theorem}

According to \cite{baierzhao}, Theorem \ref{heckesquares} and \eqref{piatbound} imply the following result on the sign changes of $\lambda(p)$ at
Piatetski-Shapiro primes $p$.

\begin{theorem} \label{hecqeoscillations} Let $1<c<25/24$ be fixed and assume that \eqref{LHLG} holds. Then 
$\lambda(p)$ changes sign infinitely often as $p$ runs through the primes of the form $p=\left[n^c\right]$ with $n\in \mathbb{N}$.
\end{theorem}

We point out that the full strength of the Lindel\"of hypothesis is not required to obtain non-trivial bounds for the exponential sums in question. However, in this paper, we want to establish the strongest possible result that our method allows.

\section{Farey dissection}
Our goal is to establish a non-trivial bound for the exponential sum
$$
\sum\limits_{n\sim N} a_n e(f(n))
$$
in Theorem 1. To this end, we shall split this exponential sum into short subsums using a Farey dissection of a certain interval. 
We note that the splitting of the summation interval in the present paper differs from that in \cite{baierzhao}.  
It will become clear in the next section why it is advantageous to split the summation interval as described below. 

For $x\ge 1$, let 
\begin{equation} \label{hdef}
h(x):=f'(x)+xf''(x).
\end{equation}
By the condition f) on $f$ in section 1, we have 
\begin{equation} \label{hprime}
h'(x)=2f''(x)+xf'''(x)\asymp \frac{f(x)}{x^2}.
\end{equation}
Hence, $h(x)$ is monotonically increasing or decreasing. In the sequel, we assume without loss of generality that $h(x)$ is 
monotonically decreasing (in particular, if $f(x)$ is defined as in \eqref{concretef} in section 10, then $h(x)$ will have this property). Let $Q$ be a real parameter with
\begin{equation} \label{Qcond}
1\le Q\le N,
\end{equation}
to be chosen later. 

Now we make a Farey dissection of level $Q$ of the interval $[h(N'),h(N))$ (for details on Farey intervals, see
\cite{baierzhao}, for example). In this way, we write $[h(N'),h(N))$ as the disjoint union of intervals of the form 
$$
\left[\left.\frac{l}{q}-\frac{M_1}{qQ},\frac{l}{q}+\frac{M_2}{qQ}\right)\right. \cap [h(N'),h(N)),
$$
where $M_1,M_2\asymp 1$, $q\le Q$ and $(q,l)=1$. Projecting these intervals back into $(N,N']$ under the map $h^{-1}$, we get intervals of the form
$$
h^{-1}\left(\left.\left[\frac{l}{q}-\frac{M_1}{qQ},\frac{l}{q}+\frac{M_2}{qQ}\right.\right) \cap [h(N'),h(N))\right)=
(x_0-m_1,x_0+m_2] \subseteq (N,N']
$$
with 
$$
x_0=h^{-1}\left(\frac{l}{q}\right)
$$
and
$$
m_1,m_2 \ll \frac{1}{qQ} \cdot \left(h^{-1}\right)'\left(\frac{l}{q}\right)=\frac{1}{qQ}\cdot \frac{1}{h'(x_0)} \asymp 
\frac{N^2}{qQf(N)},
$$
by \eqref{hprime} and the conditions b) and c) on  $f$ in section 1. 

In the following sections, we shall estimate the subsums
$$
\sum\limits_{x_0-m_1<n\le x_0+m_2} a_n e(f(n)).
$$

\section{Approximation of $f$ and partial summation}
In $(x_0-m_1,x_0+m_2]$, we now approximate the function $f(x)$ by
\begin{equation} \label{gdef}
g(x)=h(x_0)x-x_0^2f''(x_0)\log x+C=\frac{l}{q}\cdot x-x_0^2f''(x_0)\log x+C,
\end{equation}
where 
$$
C:=f(x_0)-h(x_0)x_0+x_0^2f''(x_0)\log x_0.
$$
Using the definition of $h(x)$ in \eqref{hdef}, It follows that
$$
g(x_0)=f(x_0),\quad g'(x_0)=f'(x_0), \quad \mbox{and} \quad
g''(x_0)=f''(x_0).
$$

Hence, applying Taylor's theorem to approximate $(f-g)'(x)$ near $x_0$, we have 
$$
f'(x)-g'(x)=\frac{1}{2}(x-x_0)^2\left(f'''(c)-g'''(c)\right)
$$
for some $c\in \left[x_0-m_1,x_0+m_2\right]$ if $x\in (x_0-m_1,x_0+m_2]$. Now,
$$
f'''(c)-g'''(c)=f'''(c)-\frac{2x_0^2f''(x_0)}{c^3}\ll \frac{f(N)}{N^3},
$$
by our conditions on $f$. Hence,
$$
f'(x)-g'(x)\ll \frac{N^4}{q^2Q^2f(N)^2}\cdot \frac{f(N)}{N^3}\ll \frac{N}{q^2Q^2f(N)}.
$$

Using partial summation, we deduce that
\begin{equation} \label{lange}
\begin{split}
&  \sum\limits_{n\in (x_0-m_1,x_0+m_2]} a_n e(f(n))\\
= & \sum\limits_{n\in (x_0-m_1,x_0+m_2]} a_n e(g(n))e(f(n)-g(n))\\
= & e(f(x_0+m_2)-g(x_0+m_2))\sum\limits_{n\in (x_0-m_1,x_0+m_2]} a_n e(g(n))-\\ & 2\pi i 
\int\limits_{x_0-m_1}^{x_0+m_2} \left(\sum\limits_{n\in(x_0-m_1,u]} a_ne(g(n))\right)(f'(u)-g'(u)) 
e(f(u)-g(u)) \dif u\\
\ll & \left(1+(m_1+m_2)\cdot \frac{N}{q^2Q^2f(N)}\right) \cdot \max_{u\le x_0+m_2} \left| 
\sum\limits_{n\in (x_0-m_1,u]} a_ne(g(n)) \right|\\
\ll & \left(1+\frac{N^3}{q^3Q^3f(N)^2}\right) \cdot \max_{u\le x_0+m_2} \left| \sum\limits_{n\in (x_0-m_1,u]} a_ne(g(n))
\right|.
\end{split}
\end{equation}
Thus we have replaced the function $f(x)$ by $g(x)$. The exponential sum with $g(n)$ in place of $f(n)$ can now be related to the functions
$F(s,\chi)$. This will be done in the next sections. 

In \cite{baierzhao}, we approximated the function $f(x)$ just by a linear function of the form 
$$
g(x)=\frac{l}{q}\cdot x +C
$$
in an interval around the point $x_0=f'^{-1}(l/q)$. However, in this way one can just force the first derivative of $g(x)$ 
to agree with that of $f(x)$ at the point at $x=x_0$. The approximation of $f(x)$ by the function $g(x)$ defined in \eqref{gdef} allows 
to force the first {\it and the second} derivatives of $g(x)$ and $f(x)$ to 
agree at $x=x_0$.  This reduces the error in the approximation substantially and is the key point of this paper.  

\section{Rewriting $\sum_{n} a_n e(g(n))$ using multiplicative characters}
We have
\begin{equation} \label{qurze}
\sum\limits_{x_0-m_1< n\le u} a_n e(g(n))=\sum\limits_{x_0-m_1< n\le u} a_n e\left(n\cdot \frac{l}{q}\right) \cdot
n^{-iT} 
\end{equation}
with
\begin{equation} \label{Tchoice}
T:=2\pi x_0^2 f''(x_0).
\end{equation}
We break the sum over $n$ as follows.
\begin{eqnarray*}
\sum\limits_{x_0-m_1< n\le u} a_n e\left(n\cdot \frac{l}{q}\right) \cdot n^{-iT}
&=& \sum\limits_{d|q} \sum\limits_{\substack{x_0-m_1< n\le u\\ (n,q)=d}} a_n e\left(n\cdot \frac{l}{q}\right) \cdot n^{-iT}\\
&=& \sum\limits_{d|q} d^{-iT} \sum\limits_{\substack{(x_0-m_1)/d< n\le u/d\\ (n,q/d)=1}} a_{dn} e\left(n\cdot \frac{l}{q/d}\right) 
\cdot n^{-iT}.
\end{eqnarray*}
Now we write the additive character in the last line using multiplicative characters in the form
$$
e\left(n\cdot \frac{l}{q/d}\right)=\frac{1}{\varphi(q/d)} \cdot \sum\limits_{\chi \bmod{q/d}} \overline{\chi}(l) \tau(\overline{\chi})
\chi(n).
$$
It follows that
\begin{equation} \label{splitting}
\begin{split}
& \sum\limits_{x_0-m_1\le n\le u} a_n e\left(n\cdot \frac{l}{q}\right) \cdot n^{-iT} \\
=& \sum\limits_{d|q} d^{-iT} \sum\limits_{\chi \bmod{q/d}} \frac{1}{\varphi(q/d)} \cdot \overline{\chi}(l)\tau(\overline{\chi})
\sum\limits_{(x_0-m_1)/d< n\le u/d} a_{dn}\chi(n)n^{-iT}.
\end{split}
\end{equation}

\section{Reduction to $F(s,\chi)$}
Using Perron's formula and the Ramanujan conjecture, $a_n\ll n^{\varepsilon}$, we have
\begin{equation} \label{perronsformula}
\begin{split}
\sum\limits_{(x_0-m_1)/d< n\le u/d} a_{dn}\chi(n)n^{-iT} &= \int\limits_{c-iT_0}^{c+iT_0} 
\left(\sum_{n=1}^{\infty} a_{dn}\chi(n)n^{-s-iT}\right) 
\left(\left(\frac{u}{d}\right)^{s}-\left(\frac{x_0-m_1}{d}\right)^s\right) \frac{\dif s}{s} + \\
& O\left(\frac{N^{1+\varepsilon}}{dT_0}+N^{\varepsilon}\right)
\end{split}
\end{equation}
for $c=1+1/\log N$ and $T_0\ge 1$, where we recall that 
\begin{equation} \label{intervalbounds}
N\le x_0-m_1<u\le N'\le 2N.
\end{equation}
Next, we relate the Dirichlet series in the integrand to $F(s,\chi)$.

Since $F(s)$ has an Euler product, the coefficients $a_n$ of $F(s)$ are multiplicative in $n$. Hence,
for $\Re s>1$, we have
$$
\sum\limits_{n=1}^{\infty} a_{dn}\chi(n)n^{-s} = \left(\sum\limits_{\substack{n=1\\ s(n)|d}}^{\infty} a_{dn}\chi(n)n^{-s}\right)
\cdot \left(\sum\limits_{\substack{n=1\\ (n,d)=1}}^{\infty} a_{dn}\chi(n)n^{-s}\right),
$$
where $s(n)$ is the largest squarefree number dividing $n$, and we may write
$$
\sum\limits_{\substack{n=1\\ s(n)|d}}^{\infty} a_{dn}\chi(n)n^{-s}=\prod\limits_{p|d} \sum\limits_{k=0}^{\infty} 
a_{p^{\alpha(p)+k}}\chi^k(p)p^{-ks},
$$ 
where 
$$
d=\prod\limits_{p|d} p^{\alpha(p)}
$$
is the prime number factorization of $d$. Further,
$$
\sum\limits_{\substack{n=1\\ (n,d)=1}}^{\infty} a_n\chi(n)n^{-s}=
F(s,\chi)\prod\limits_{p|d} \left(\sum\limits_{k=0}^{\infty} 
a_{p^k}\chi^k(p)p^{-ks}\right)^{-1}.
$$
So altogether,
$$
\sum\limits_{n=1}^{\infty} a_{dn}\chi(n)n^{-s}=G_d(s,\chi)F(s,\chi)
$$
with
$$
G_d(s,\chi)=\prod\limits_{p|d} \frac{\sum\limits_{k=0}^{\infty} 
a_{p^{\alpha(p)+k}}\chi^k(p)p^{-ks}}{\sum\limits_{k=0}^{\infty} 
a_{p^k}\chi^k(p)p^{-ks}}.
$$

Hence, the integral on the right-hand side of \eqref{perronsformula} takes the form
\begin{equation} \label{inttranslate}
\begin{split}
& \int\limits_{c-iT_0}^{c+iT_0} \left(\sum_{n=1}^{\infty} a_{dn}\chi(n)n^{-s-iT}\right) 
\left(\left(\frac{u}{d}\right)^{s}-\left(\frac{x_0-m_1}{d}\right)^s\right) \frac{\dif s}{s}\\ &=
\int\limits_{c-iT_0}^{c+iT_0} G_d(s+iT,\chi)F(s+iT,\chi)
\left(\left(\frac{u}{d}\right)^{s}-\left(\frac{x_0-m_1}{d}\right)^{s}\right) \frac{\dif s}{s}.
\end{split}
\end{equation}

\section{Estimation of the integral}
We shall need a bound for $G_d(s,\chi)$ if $\Re s\ge 1/2$, which we establish in the following. 
By the Ramanujan conjecture, $a_n\ll n^{\varepsilon}$, we have
$$
\sum\limits_{k=0}^{\infty} a_{p^{\alpha(p)+k}}\chi^k(p)p^{-ks} \ll p^{\alpha(p)\varepsilon} \quad \mbox{uniformly for } 
\Re s\ge \frac{1}{2}. 
$$

Let $t_0\in \mathbb{R}$ such that 
$$
p^{-it_0}=\chi(p).
$$
Then
$$
\sum\limits_{k=0}^{\infty} a_{p^k}\chi^k(p)p^{-ks}=\sum\limits_{k=0}^{\infty} a_{p^k}p^{-k(s+it_0)},
$$
and by axiom (iii) (Euler product) for $F(s)$, we therefore have
$$
\log \sum\limits_{k=0}^{\infty} a_{p^k}\chi^k(p)p^{-ks} = \log F_p(s+it_0) = \sum_{k=0}^\infty b_{p^k}p^{-k(s+it_0)}=
\sum_{k=0}^\infty b_{p^k}\chi^k(p)p^{-ks},
$$
where $b_{p^k}$ ($k=0,1,...$) are suitable coefficients satisfying 
$$
b_{p^k}\ll p^{k\theta}
$$
for some $\theta<1/2$. It follows that
\begin{equation} \label{-1}
\left(\sum\limits_{k=0}^{\infty} a_{p^k}\chi^k(p)p^{-ks}\right)^{-1}=O(1) \quad \mbox{uniformly for $\Re s\ge 1/2$.}
\end{equation}

From the above estimates, we deduce that
$$
G_d(s,\chi)\ll d^{\varepsilon} \quad \mbox{uniformly for $\Re s\ge 1/2$.}
$$
Furthermore, using conditions b) and c) on $F(s,\chi)$ in section 1 together with the Phragmen-Lindel\"of principle, we have
$$
F(\sigma+it,\chi)\ll (tq)^{\varepsilon} \quad \mbox{uniformly for } \sigma\ge \frac{1}{2} \mbox{ and } |\sigma+it-1|\ge \frac{1}{4}.
$$
Thus,
\begin{equation} \label{GFestimate}
G_d(\sigma+it,\chi)F(\sigma+it,\chi)\ll (dtq)^{\varepsilon} \quad \mbox{uniformly for } \sigma\ge \frac{1}{2} \mbox{ and } 
|\sigma+it-1|\ge \frac{1}{4}.
\end{equation}

Now we bound the integral on the right-hand side of \eqref{inttranslate}, where we suppose that 
\begin{equation} \label{Tcond}
T_0\le \frac{T}{2}.
\end{equation}
Using Cauchy's integral theorem, we then have
\begin{equation} \label{Cauchy}
\begin{split}
& \int\limits_{c-iT_0}^{c+iT_0} G_d(s+iT,\chi)F(s+iT,\chi)
\left(\left(\frac{u}{d}\right)^{s}-\left(\frac{x_0-m_1}{d}\right)^{s}\right) \frac{\dif s}{s}\\ = &
\left(\int\limits_{c-iT_0}^{1/2-iT_0}+\int\limits_{1/2-iT_0}^{1/2+iT_0}+
\int\limits_{1/2+iT_0}^{c+iT_0}\right) G_d(s+iT,\chi)F(s+iT,\chi)\times\\ &
\left(\left(\frac{u}{d}\right)^{s}-\left(\frac{x_0-m_1}{d}\right)^{s}\right) \frac{\dif s}{s}.
\end{split}
\end{equation}
From \eqref{intervalbounds}, \eqref{GFestimate}, \eqref{Tcond}, \eqref{Cauchy} and $d|q$, we deduce that
\begin{equation} \label{intest}
\int\limits_{c-iT_0}^{c+iT_0} G_d(s+iT,\chi)F(s+iT,\chi)
\left(\left(\frac{u}{d}\right)^{s}-\left(\frac{x_0-m_1}{d}\right)^{s}\right) \frac{\dif s}{s}
\ll (Tq)^{\varepsilon}\cdot \left(\left(\frac{N}{d}\right)^{1/2} + \frac{N}{dT_0}\right)
\end{equation}
if $T\ge 1/2$.

\section{Proof of Theorem 1}
Now we choose 
\begin{equation} \label{T0choice}
T_0:=\left(\frac{N}{d}\right)^{1/2}.
\end{equation}
We note that
\begin{equation} \label{Tasymp}
T\asymp f(N)
\end{equation}
by \eqref{Tchoice}, $x_0\asymp N$ and our conditions of $f$.
Hence, by \eqref{anotherfcond}, the condition \eqref{Tcond} is satisfied if $N$ is large enough. From \eqref{anotherfcond}, \eqref{perronsformula}, \eqref{inttranslate}, \eqref{intest}, \eqref{T0choice}, \eqref{Tasymp} and $d\le N$ (by $d|q$, $q\le Q$ and \eqref{Qcond}), we deduce that
$$
\sum\limits_{(x_0-m_1)/d< n\le u/d} a_{dn}\chi(n)n^{-iT} \ll (qN)^{\varepsilon} \cdot \left(\frac{N}{d}\right)^{1/2}.
$$
Plugging this into \eqref{splitting}, and using $\varphi(q/d)\gg q^{1-\varepsilon}/d$ and $|\tau(\overline{\chi})|\le \sqrt{q/d}$, we obtain
\begin{equation}
\left|\sum\limits_{x_0-m_1< n\le u} a_n e\left(n\cdot \frac{l}{q}\right)n^{-iT}\right| \ll (qN)^{1/2+\varepsilon}.
\end{equation}
This together with \eqref{lange} and \eqref{qurze} yields
\begin{equation} \label{conti}
\sum\limits_{n\in (x_0-m_1,x_0+m_2]} a_n e(f(n))\ll \left(1+\frac{N^3}{q^3Q^3f(N)^2}\right)(qN)^{1/2+\varepsilon}.
\end{equation}

In section 3, we have divided the interval $[h(N'),h(N))$ into Farey intervals around fractions $l/q$ with
$$
1\le q\le Q, \quad l\asymp q\cdot h(N) \asymp q\cdot \frac{f(N)}{N} \quad \mbox{and} \quad (q,l)=1.  
$$
Hence, summing the contributions of the short sums in \eqref{conti} over all relevant $q$ and $l$, we get
\begin{eqnarray*}
\sum\limits_{n\sim N} a_n e(f(n))
&\ll& \sum\limits_{q\le Q}\ \sum\limits_{l\asymp qf(N)/N} \left(1+\frac{N^3}{q^3Q^3f(N)^2}\right)(qN)^{1/2+\varepsilon}\\
&\ll& \left(\frac{Q^{5/2}f(N)}{N^{1/2}}+\frac{N^{5/2}}{Q^3f(N)}\right)(QN)^{\varepsilon}. 
\end{eqnarray*}
Now we choose
$$
Q:=\left(\frac{N^{5/2}}{f(N)}\cdot \frac{N^{1/2}}{f(N)}\right)^{2/11}=\frac{N^{6/11}}{f(N)^{4/11}}.
$$
Thus we get
$$
\sum\limits_{n\sim N} a_ne(f(n)) \ll N^{19/22+\varepsilon}f(N)^{1/11},
$$
which completes the proof. $\Box$

\section{Proof of Theorem 2}
Theorem 2 can be proved along similar lines as Theorem 1.
The arguments in sections 3-5 carry over completely. We are then led to the sum
$$
\sum\limits_{\substack{(x_0-m_1)/d< n\le t/d\\ (n,m)=1}} \mu^2(dn)a_{dn} \chi(n)n^{-iT}
$$
in place of the sum
$$
\sum\limits_{(x_0-m_1)/d< n\le t/d} a_{dn} \chi(n)n^{-iT}
$$
considered in sections 6-8.
We now use the fact that $a_{n}$ is multiplicative in $n$ to rewrite the sum in question in the form
$$
\sum\limits_{\substack{(x_0-m_1)/d< n\le t/d\\ (n,m)=1}} \mu^2(dn)a_{dn} \chi(n)n^{-iT}=\mu^2(d) a_d 
\sum\limits_{(x_0-m_1)/d< n\le t/d} \mu^2(n)a_{n} \chi_1(n)n^{-iT},
$$
where $\chi_1(n)=\chi(n)\chi_0(n)$, $\chi_0(n)$ being the principal character modulo $dm$. Similarly as in section 6,
we relate the sum over $n$ on the right-hand side to the corresponding Dirichlet series, which we write in the form
$$
\sum\limits_{n=1}^{\infty} \mu^2(n)a_n\chi_1(n)n^{-s}=\prod\limits_{p}\left(1+a_p\chi_1(p)p^{-s}\right)=H(s,\chi_1)F(s,\chi_1),
$$
where
$$
H(s,\chi_1)=\prod\limits_{p} \frac{1+a_p\chi_1(p)p^{-s}}{\sum\limits_{k=0}^{\infty} a_{p^k}\chi_1^k(p)p^{-ks}}=
\prod\limits_{p} \left(1-\frac{\sum\limits_{k=2}^{\infty} a_{p^k}\chi_1^k(p)p^{-ks}}{\sum\limits_{k=0}^{\infty} a_{p^k}\chi_1^k(p)p^{-ks}}\right).
$$
By $a_n\ll n^{\varepsilon}$ and \eqref{-1}, the product on the right-hand side 
converges absolutely and uniformly in every compact subset $S$ 
of the half plane $\Re s>1/2$. Hence, the function $H(s,\chi_1)$ is entire there.
Moreover, $|H(s,\chi_1)|$ is bounded by a constant $C(\varepsilon)$ if $\Re s\ge 1/2+\varepsilon$. The rest of the proof
follows the arguments in the proof of Theorem 1, where the function $G_d(s,\chi)$ is replaced by $H(s,\chi_1)$, and in the
application of Cauchy's integral theorem, the line of integration is shifted to $\Re s=1/2+\varepsilon$ instead of $\Re s=1/2$.

\section{Proof of Theorem 3}
The general procedure of the proof will be similar as in \cite{baierzhao}, where we bounded the sum
$$
\sum\limits_{\substack{n\le N\\ \left[n^c\right]\in \mathbb{P}}} \lambda\left(\left[n^c\right]\right).
$$
Therefore, we will be very brief in general and go into details only in the parts where the proof of Theorem 3 deviates substantially from
that of Theorem 1 in \cite{baierzhao}. First, we use the well-known relation
$$
\lambda(p)^2=1+\lambda\left(p^2\right).
$$
Hence, we have
\begin{equation} \label{sums}
\sum\limits_{\substack{n\le N\\ \left[n^c\right]\in \mathbb{P}}} \lambda\left(\left[n^c\right]\right)^2 = 
\sum\limits_{\substack{n\le N\\ \left[n^c\right]\in \mathbb{P}}} 1 + \sum\limits_{\substack{n\le N\\ \left[n^c\right]\in \mathbb{P}}} \lambda\left(\left[n^c\right]^2\right).
\end{equation}
The ordinary Piatetski-Shapiro prime number theorem (see \cite{baierzhao}, for example) tells us that
$$
\sum\limits_{\substack{n\le N\\ \left[n^c\right]\in \mathbb{P}}} 1 \sim \frac{N}{c\log N} \quad \mbox{as } N\rightarrow \infty
$$
for every fixed $c$ in the range in Theorem 3. It remains to estimate the second sum on the right-hand side of \eqref{sums}. We write this sum in the form
$$
\sum\limits_{\substack{n\le N\\ \left[n^c\right]\in \mathbb{P}}} \lambda\left(\left[n^c\right]^2\right)=\sum\limits_{\substack{n\le N\\ \left[n^c\right]\in \mathbb{P}}} b_{\left[n^c\right]},
$$
where 
\begin{equation} \label{asin}
b_n:=\mu^2(n)\lambda\left(n^2\right)=\mu^2(n) a_n,
\end{equation}
with $a_n$ as in \eqref{theas}.

Clearly, it now suffices to bound the sum
$$
\sum\limits_{n\le N} b_{\left[n^c\right]}\Lambda\left(\left[n^c\right]\right),
$$
where $\Lambda(n)$ is the von Mangoldt function. Similarly as in \cite{baierzhao}, we pull out a main term which we estimate using an analogue of the prime number theorem for $\lambda(p^2)$. Then we reduce the error term to exponential sums as in \cite{baierzhao} and treat the von Mangoldt function appearing in them using a Vaughan-type identity due to Heath-Brown (Lemma 4 in \cite{baierzhao}). This leads to type I and type II sums.
In \cite{baierzhao}, we then used the decomposition in \eqref{multprop}
to separate the summation variables $m$ and $n$ in the said type I and type II sums.
The decomposition of $b_{mn}$ needed here is simpler since we have
$$
b_{mn}=\begin{cases} 0 & \mbox{ if } (m,n)>1,\\ b_mb_n & \mbox{ if } (m,n)=1\end{cases}
$$
due to the appearance of the M\"obius function in the definition of $b_n$. Now the type I and type II sums take the form
\begin{equation} \label{K}
K=\sum\limits_{h\sim H} \mathop{\sum\limits_{m\sim X}\sum\limits_{n\sim Y}}_{(m,n)=1} C_hA_mb_ne\left(h(mn)^{\gamma}\right)
\end{equation}
and
$$
L=\sum\limits_{h\sim H} \mathop{\sum\limits_{m\sim X}\sum\limits_{n\sim Y}}_{(m,n)=1} C_hA_mB_ne\left(h(mn)^{\gamma}\right),
$$
where 
\begin{equation} \label{HXY}
\gamma=\frac{1}{c},\quad 1\le H\le N^{1-\gamma+\eta}\quad  \mbox{and}\quad  XY=N.
\end{equation} 
Here $C_h$, $A_m$ and $B_n$ are general coefficients of size $\ll N^{\varepsilon}$, and $b_n$ is defined as in  \eqref{asin}.

We remove the coprimality condition $(m,n)=1$ in $K$ and $L$ using M\"obius inversion, getting
\begin{equation} \label{L}
K=\sum\limits_{d} \mu(d)K_d \quad \mbox{and} \quad L=\sum\limits_{d} \mu(d)L_d
\end{equation}
with
$$
K_d=\sum\limits_{h\sim H} \ \sum\limits_{m\sim X/d} \ \sum\limits_{n\sim Y/d} C_hA_{dm}b_{dn}e\left(h\left(d^2mn\right)^{\gamma}\right).
$$
and 
$$
L_d=\sum\limits_{h\sim H} \ \sum\limits_{m\sim X/d} \ \sum\limits_{n\sim Y/d} C_hA_{dm}B_{dn}e\left(h\left(d^2mn\right)^{\gamma}\right).
$$
Using \eqref{L} and Lemmas 15, 16 and 18 in \cite{baierzhao}, we deduce that
$$
L\ll N^{1-\eta} \quad \mbox{if } N^{1-\gamma+100\eta}\le Y\le N^{5\gamma-4-100\eta}
$$
and
\begin{equation} \label{K1}
K\ll N^{1-\eta} \quad \mbox{if } N^{3-3\gamma+100\eta}\le Y\le N^{\gamma-100\eta},
\end{equation}
for some small $\eta>0$, provided that $\gamma>7/8$.

We note that if $0<\gamma<1$, then the function
\begin{equation} \label{concretef}
f(x)=h(mx)^{\gamma}
\end{equation}
satisfies the conditions on $f$ in section 1, and if $\gamma>1/2$, $h\sim H$, $m\sim X$, $Y\ge N^{2/3+100\eta}$, $\eta$ is sufficiently small and \eqref{HXY} is satisfied, then condition
\eqref{anotherfcond}, with $N$ replaced by $Y$, in Theorems 1 and 2 holds, i.e.
$$
Y^{1/2+\eta}\le f(Y)\le Y^{3/2-\eta}.
$$
Now applying Theorem 2 with $f(n)$ defined as above and $a_n$ defined as in \eqref{theas} to the inner sum over $n$ on the right-hand side of \eqref{K}, we get
$$
\sum\limits_{\substack{n\sim Y\\ (m,n)=1}} b_ne\left(h(mn)^{\gamma}\right)=\sum\limits_{\substack{n\sim Y\\ (m,n)=1}} \mu^2(n)a_ne\left(h(mn)^{\gamma}\right)\ll
H^{1/11}X^{\gamma/11}Y^{19/22+\gamma/11+\varepsilon}
$$
and hence
$$
K\ll H^{12/11}X^{1+\gamma/11}Y^{19/22+\gamma/11+\varepsilon}, 
$$
provided that $1/2<\gamma<1$ and $Y\ge N^{2/3+100\eta}$.
By \eqref{HXY}, it follows that 
\begin{equation} \label{K2}
K\ll N^{1-\eta} \quad \mbox{if } Y\ge N^{8-22\gamma/3+100\eta},
\end{equation}
provided that $1/2<\gamma<1$.
The $Y$-ranges in \eqref{K1} and \eqref{K2} overlap if $\gamma> 24/25$
and $\eta$ is small enough. Hence, we have
$$
K\ll N^{1-\eta} \quad \mbox{if } Y\ge N^{3-3\gamma+100\eta}, 
$$
provided that $24/25<\gamma<1$ and $\eta$ is sufficiently small. The rest of the proof is similar as in section 12 in \cite{baierzhao}. We note that the range $1<c<25/24$ in Theorem 3
comes from the above condition $24/25<\gamma<1$.\\

\noindent{\bf Acknowledgments.} This work was supported by an ERC grant 258713.

$ $\\
\noindent\begin{tabular}{p{8cm}p{8cm}}
Stephan Baier\\
Mathematisches Institut\\
Universit\"at G\"ottingen\\
Bunsenstr.\ 3--5, 37073\\
G\"ottingen Germany\\
Email: {\tt sbaier@uni-math.gwdg.de}
\end{tabular}
\end{document}